%% file: root.tex
\documentclass[letterpaper,10pt, journal,twoside]{IEEEtran}  

\IEEEoverridecommandlockouts                              

\usepackage{amsmath,amssymb,color,cite}
\usepackage{graphicx}

\usepackage{latexsym}
\usepackage{verbatim}
\usepackage{cite}
\usepackage[english]{babel}
\usepackage[latin1]{inputenc}
\usepackage{enumerate}
\usepackage{amsmath,amsthm}

\usepackage{tikz}
\usetikzlibrary{tikzmark}
\usepackage{amsfonts}
\usepackage{amssymb}
\usepackage{amsbsy}
\usepackage{bm}

    \let\leq\leqslant
    \let\geq\geqslant
    \let\emptyset\varnothing

\newcommand{\rrp}{\calR_+}

\input{ka-newcommands}

\newcounter{todocounter}
\usepackage[colorinlistoftodos]{todonotes}

\newtheorem{theorem}{Theorem}
\newtheorem{lemma}[theorem]{Lemma}

\theoremstyle{definition}
\newtheorem{definition}[theorem]{Definition}
\newtheorem{example}[theorem]{Example}

\newtheorem{problem}{Problem}
\newtheorem{conjecture}[theorem]{Conjecture}
\newtheorem{remark}[theorem]{Remark}

\title{Data informativity for analysis of linear systems with convex conic constraints}

\author{Jaap Eising, M. Kanat Camlibel,~\IEEEmembership{Member,~IEEE}}

\begin{document}

\maketitle

\renewcommand{\thefootnote}{\fnsymbol{footnote}}

\footnotetext{The authors are with the Jan C. Willems Center for Systems and Control and the  Bernoulli Institute for Mathematics, Computer Science, and Artificial Intelligence, University of Groningen, Nijenborgh 9, 9747 AG, Groningen, The Netherlands. (email: {\tt j.eising@rug.nl; m.k.camlibel@rug.nl})}

\begin{abstract}
	This paper studies the informativity problem for reachability and null-controllability of constrained systems. To be precise, we will focus on an unknown linear systems with convex conic constraints from which we measure data consisting of exact state trajectories of finite length. We are interested in performing system analysis of such an unknown system on the basis of the measured data. However, from such measurements it is only possible to obtain a unique system explaining the data in very restrictive cases. This means that we can not approach this problem using system identification combined with model based analysis. As such, we will formulate conditions on the data under which any such system consistent with the measurements is guaranteed to be reachable or null-controllable. These conditions are stated in terms of spectral conditions and subspace inclusions, and therefore they are easy to verify.
	\end{abstract}
\vspace{-1em}
\section{Introduction}
This paper deals with the question: what can be inferred from an unknown constrained linear system on the basis of state measurements? A similar question, for unconstrained systems, has recently led to the development of the \textit{informativity framework} in \cite{vanWaarde2020}. The observation at the center of this framework is that we can only conclude that the unknown system has a given property if \textit{all} systems compatible with the measurements have this property. In the context of linear systems this has lead to, among others, results for analysis problems in \cite{Eising2020c} and control problems in \cite{vanWaarde2020c,Trentelman2020}. Parallel to the work performed within this framework, similar analysis problems are addressed in \cite{Koch2020}, while control problems are addressed in \cite{DePersis2020,Berberich2020}. 
            

In contrast to this earlier work, we will be focusing on conically constrained linear systems. Such conic constraints often arise naturally in modeling, taking the form of e.g. nonnegativity constraints on the input or states. Specifically, we will be looking at the class of difference inclusions of the form
\[ x_{k+1} \in H(x_k)\]
where $H:\mathbb{R}^n\rightrightarrows \mathbb{R}^n$ is a convex process, that is, a set-valued map whose graph is a convex cone. It it straightforward to show that any conically constrained linear system can be written as such a system and vice versa. Such difference inclusions arise naturally in many different contexts, including chemical reaction networks \cite{Angeli:09}, von Neumann-Gale economic growth models \cite{makarov:77} and cable-suspended robots \cite{j37,oh:05-2}. Lastly, as shown in e.g. \cite{Frankowska:87,Aubin:84}, difference inclusions of convex processes can be used as meaningful approximations of more complex set-valued maps. 

The many applications of convex processes have led to interest in the analysis of such systems. In particular, this paper will consider the system-theoretic properties of reachability and null-controllability. For a given convex process, tests for these properties have been developed in terms of spectral conditions. First among these were the characterizations of reachability and null-controllability in \cite{AFO:86,PD:94}. However, the aforementioned characterizations only regard \textit{strict} (nonempty everywhere) convex processes, which limits the applicability for our goals. In \cite{ReachNullc:19} both of these results are generalized to work for a class of nonstrict convex processes. These characterizations of \cite{ReachNullc:19} will be fundamental in our investigation of informativity. In this paper, we will be interested in analysing whether these system-theoretic properties hold for \textit{all} convex processes compatible with measured data. Specifically, the data we consider will consist of exact measurements of the state. This means that we leave studies involving measurement noise or partial state measurements as future extensions.


Apart from the aforementioned work, some results in data-driven analysis and control should be mentioned. With regard to unconstrained linear systems \cite{Park2009} analyses stability of an input/output system using time series data. The works \cite{Wang2011,Liu2014,Niu2017,Zhou2018} deal with data-based controllability and observability analysis. Lastly, many methods arising from Model Predictive Control (MPC) are well suited to constrained systems. For an overview of such methods, we refer to \cite{Mayne:00,Mayne2014}. More specifically, MPC has recently been brought into a data-based context in \cite{Piga2018,Coulson2019}.

The contribution of this paper is threefold: 
\begin{enumerate}
	\item We expand the informativity framework of \cite{vanWaarde2020} towards the class of convex processes. This framework will naturally lead to the formulation of a number of problems. In particular, we will illustrate the framework by resolving the problems of informativity for reachability and null-controllability.
	\item We develop explicit tools to manipulate and perform analysis on convex processes with a polyhedral graph. Assuming polyhedrality will allow us to represent convex processes and the conditions required for reachability and null-controllability in a convenient way. 
	\item Lastly, we note the fact that polyhedral convex processes naturally arise from the aforementioned informativity problems with finite measurements. This allows us to combine the previous points to formulate tests \textit{on measured state data} to conclude that all convex processes consistent with the data are reachable or null-controllable. 
\end{enumerate} 

This paper is organized as follows: We begin in Section~\ref{sec:cp} with definitions of convex process and reachability and null-controllability. After this, Section~\ref{sec:prob} introduces informativity and formally states the problem we will consider in this paper. In Section~\ref{sec:prelim cp}, we will present some known results regarding the analysis of convex processes, which will be applied in Section~\ref{sec:inf} to our problem. We finalize the paper with conclusions in Section~\ref{sec:conc}.
\vspace{-0.5em}
\section{Convex processes}\label{sec:cp} 

Given convex sets $\calS,\calT\subseteq\mathbb{R}^n$ and scalar $\rho\in\R$ we define the sum and scalar multiplication of sets as: 
\[ \calS+\calT := \{s+t \mid s\in\calS, t\in\calT\}, \quad \rho \calS := \{\rho s \mid s\in\calS\}. \] 
We denote the closure of $\calS$ by $\cl\calS$. A \textit{convex cone} is a nonempty convex set that is closed under nonnegative scalar multiplication. 

A \textit{set-valued map}, denoted $H:\R^n\rightrightarrows \R^n$ is a map taking elements of $\mathbb{R}^n$ to subsets of $\mathbb{R}^n$. It is called a {\em convex process\/}, {\em closed convex process\/} or {\em linear process\/} if its graph
\[ \graph H := \{(x,y)\in \R^n\times \R^n \mid y\in H(x)\}\]
is a convex cone, closed convex cone or subspace, respectively.

The \textit{domain} and \textit{image} of $H$ are defined as $\dom H= \{x\in \R^n\mid H(x)\neq \emptyset \}$ and $\im H = \{y\in \R^n : \exists\, x \textrm{ s.t. } y\in H(x)\}$. If $\dom H=\R^n$, we say that $H$ is \textit{strict}.

In this paper we consider systems described by a \textit{difference inclusion} of the form
\vspace{-0.5em}\begin{equation}\label{eq:diffinc} x_{k+1} \in H(x_k)
\end{equation}
where $H:\R^n\rightrightarrows \R^n$ is a convex process. Our main motivation for considering this class of systems is the fact that this class of systems captures the behavior of all linear systems with convex conic constraints. This will be made explicit in the following example. 

\begin{example}\label{ex:linear cons} 
	Consider states $x_k$ in $\mathbb{R}^n$ and inputs $u_k\in \mathbb{R}^m$. Let $A$ and $B$ be linear maps of appropriate dimensions and let $\calC\subseteq\mathbb{R}^{n+m}$ be a convex cone. Consider the linear system with conic constraints given by:
	\begin{equation}\label{eq:cons system} 
	x_{k+1}= Ax_k +Bu_k, \quad\begin{bmatrix} x_k \\ u_k \end{bmatrix} \in \calC. 
	\end{equation}
	Note that this description can be applied to any combination of input, state and output constraints. 
	
	We can describe the dynamics of \eqref{eq:cons system} by the difference inclusion~\eqref{eq:diffinc} with the convex process $H$ defined by:
	\[  H(x) := \left\lbrace Ax+Bu \,\middle\vert\,  \begin{bmatrix} x \\ u \end{bmatrix} \in \calC \right\rbrace. \]
	This reveals that we can study the properties of conically constrained linear systems by studying convex processes, without any loss of generality.
\end{example}

Next, we define a number of sets associated with the difference inclusion \eqref{eq:diffinc}. A $q$\textit{-step trajectory} is a (finite) sequence $x_0,\ldots,x_q$ such that \eqref{eq:diffinc} holds for all $k<q$. We define the $q$\textit{-step behavior} as:
\[ \mathfrak{B}_q(H) := \left\lbrace (x_k)_{k=0}^q \in (\mathbb{R}^n)^{q+1}\mid (x_k) \textrm{ satisfies } (\ref{eq:diffinc})\right\rbrace . \] 
Using this, we define the \textit{reachable} and \textit{null-controllable} sets by:
\bse\label{eq:reach null def}
\begin{align*}
\!\!\calR(H) &\!:=\! \big\lbrace \xi \mid \exists q, (x_k)_{k=0}^q \in\mathfrak{B}_q(H) \textrm{ s.t. } x_0=0, x_q=\xi \big\rbrace, \\ 
\!\!\!\calN(H) &\!:=\! \big\lbrace \xi \mid \exists q, (x_k)_{k=0}^q \in\mathfrak{B}_q(H)\textrm{ s.t. } x_0=\xi, x_q=0 \big\rbrace.
\end{align*}
\ese
We say that a point $\xi\in\mathbb{R}^n$ is \textit{reachable} if $\xi\in\calR(H)$. That is, there exists a $q$-step trajectory from the origin to $\xi$. Similarly, we say a point $\xi\in\mathbb{R}^n$ is \textit{null-controllable} if $\xi\in\calN(H)$. 

By a \textit{trajectory} of \eqref{eq:diffinc}, we mean a sequence $(x_k)_{k\in\N }$ such that \eqref{eq:diffinc} holds for all $k\geq 0$. The \textit{behavior} is the set of all trajectories: 
\[ \mathfrak{B}(H) := \left\lbrace (x_k) \in (\mathbb{R}^n)^\mathbb{N}\mid (x_k) \textrm{ is a trajectory of } \eqref{eq:diffinc}\right\rbrace . \]

The set of \textit{feasible} states of the difference inclusion \eqref{eq:diffinc} is the set of states from which a trajectory emanates:	
\[ \mathcal{F}(H) := \left\lbrace\xi \mid \exists (x_k)\in \mathfrak{B}(H) \textrm{ with } x_0=\xi \right\rbrace.\] 
Clearly, if $H$ is a convex process, then $\calF(H)$ is a convex cone.

It is important to stress that in general not every point in the state space is feasible: In Example~\ref{ex:linear cons}, if we consider a point $x_0$ for which no $u_0$ satisfies the constraints, we have that $H(x_0)=\emptyset$. This means that $x_0$ is not a feasbile point. As there is no need to reach or control states that violate the constraints, we say the system \eqref{eq:diffinc} is \textit{reachable} or \textit{null-controllable} if every feasible state is reachable or null-controllable respectively. In terms of the previously defined sets, these can be written as $\calF(H)\subseteq\calR(H)$ and $\calF(H)\subseteq\calN(H)$ respectively. 	

It is important to note that, as is the case for discrete-time linear systems, reachability and null-controllability are not equivalent notions. 
\vspace{-0.5em}\section{Problem formulation}\label{sec:prob}
In this paper we are interested in analyzing the properties of an unknown system based on measurements performed on it. We will assume that the system under consideration is given by
\vspace{-0.5em}\[ x_{k+1} \in H_s(x_k)\]
where $H_s$ is an unknown convex process. However, we do have access to a number of exact state measurements corresponding to ($q-$step) trajectories of $H_s$. It is clear to see that we can view a single $q-$step trajectory as $q$ separate $1-$step trajectories. Therefore, without loss of generality, we assume that we measure single steps. That is, we are given a finite number of pairs $(x_k,y_k)\in \graph H_s$, with $k=0,\ldots, T$

Suppose that we are interested in characterizing reachability of $H_s$. As $H_s$ is unknown, it is indistinguishable from all other convex processes that could have generated the measurements. Therefore, we may only conclude that $H_s$ is reachable if \textit{all} convex processes that are compatible with the data are reachable. This motivates the following definition. Let $\Sigma$ denote the set of all convex processes $H:\mathbb{R}^n\rightrightarrows\mathbb{R}^n$ and let $\calD\subseteq\mathbb{R}^n\times \mathbb{R}^n$ be a finite set of measurements. Define the set of all convex processes compatible with these measurements by: 
\begin{equation}\label{eq:sigma D} \Sigma_\calD := \{ H\in\Sigma \mid \calD\subseteq\graph H\}. \end{equation}
Recall that, in order to characterize whether $H_s$ is reachable, we require all convex processes compatible with the measurements to be reachable. As such, we say that the data $\calD$ are \textit{informative for reachability} if every $H\in\Sigma_\calD$ is reachable. In a similar way we define \textit{informativity for null-controllability}.

Note that informativity is fundamentally a property of the data and the system class, but \textit{not} of the system $H_s$. This leads to the following problem formulation: 

\begin{problem}\label{prob: prob}
	Provide necessary and sufficient conditions on the data $\calD$ under which the data are informative for reachability or null-controllability. 
\end{problem}
\begin{remark}
	Following Example~\ref{ex:linear cons}, it is clear that all convex processes consistent with the data are reachable if and only if all conically constrained linear systems consistent with the data are reachable. As these problems are equivalent we will only focus on formulations in terms of convex processes in the remainder of this paper. 
\end{remark} 
It should be noted that, in certain cases, the informativity problem can be resolved trivially, as shown by the following example. 
\begin{example} 
	Let $n=1$, and assume that we measure the 2-step trajectory given by $x_0=0,$ $x_1=1,$ and $x_2=-1$. Then we have $\calD = \left\lbrace (0,1), (1,-1)\right\rbrace$. 
	
	Note that nonnegative scalar multiples of these measurements are also (finite step) trajectories of any convex process in $\Sigma_\calD$. As such, it is clear that for any $\alpha,\beta\geq 0$ we have 2-step trajectories $y_0=0$, $y_1=\alpha$, $y_2=-\alpha$ and $z_0=0$, $z_1=0$, $z_2=\beta$. Furthermore, the sum of two such 2-step trajectories is one as well. Therefore $(0,\alpha,\beta-\alpha)\in\mathfrak{B}_2(H)$ for any $H$ consistent with the data. As such, $\calR(H)=\mathbb{R}$ for any $H\in\Sigma_\calD$. 
\end{example}
In general, however, resolving the problem is not this straightforward. To be precise, it is made difficult by two things. First of all, apart from trivial examples, the set $\Sigma_\calD$ contains infinitely many convex processes. As such, it is usually not possible to take an approach based on identification. In addition, there may \textit{not} exist $q$ for a convex process $H$ such that
\[\!\!\calR(H) \!:=\! \big\lbrace \xi \mid \exists (x_k)_{k=0}^q \in\mathfrak{B}_q(H) \textrm{ s.t. } x_0=0, x_q=\xi \big\rbrace.\]
Therefore, testing whether a given convex process is reachable or null-controllable is a nontrivial problem in itself (see e.g. \cite{AFO:86,ReachNullc:19}).  
\vspace{-1em}\section{Analysis of convex processes}\label{sec:prelim cp} 
By definition a convex cone $\calC$ is closed under \textit{conic combinations}: If $c_1,...,c_\ell\in\calC$ then 
\[ \sum_{i=1}^\ell \alpha_ic_i \in \calC \quad \forall \alpha_i\geq 0. \]
The set of all conic combinations of a set $\calS$ is called the \textit{conic hull} and is denoted by $\cone  \calS $. If there exists a finite set $\calS$ such that $\calC = \cone\calS$ we say that $\calC$ is \textit{finitely generated} or \textit{polyhedral}. We denote the set of vectors of length $\ell$ with nonnegative and nonpositive elements by $\mathbb{R}_+^\ell$ and $\mathbb{R}_-^\ell$ respectively. Then, if $M\in\mathbb{R}^{k\times \ell}$ and $\calS$ is the set of columns of $M$, we have that 
\vspace{-0.5em}\begin{equation}\label{eq:conic hull rep} \cone \calS = M\mathbb{R}^\ell_+.\end{equation} 

For a nonempty set $\calC\subseteq \R^n $, we define the \textit{negative} and \textit{positive polar cone}, respectively,
\begin{align*}
\calC^- &:= \{ y\in\R^n \mid \langle x,y\rangle\leq 0 \hspace{1em} \forall x\in \calC\}, \\
\calC^+ &:= \{ y\in\R^n \mid \langle x,y\rangle\geq 0 \hspace{1em} \forall x\in \calC\}.
\end{align*}
Given sets $\calC$ and $\calS$, we have that $(\calC^-)^-=\cl(\cone\calC)$, and:
\begin{equation}\label{eq:dual sum} (\calC+\calS)^-  =\calC^-\cap \calS^-, \quad (\calC\cap\calS)^- = \cl (\calC^-+\calS^-). \end{equation}
Let $A$ be a linear map and let $\bullet^{-1}$ denotes the inverse image, that is, $A^{-1}\calC^- = \{ x \mid Ax\in \calC^-\}$. Then if $\calC$ is a convex cone  we have that (see e.g. \cite[Theorem 2.4.3]{Aubin:90}):
\begin{equation}\label{eq:dual and linear map}(A^\top\calC)^- = A^{-1} \calC^-.\end{equation}
The aforementioned properties also hold for the positive polar cone.

Let $H:\mathbb{R}^n\rightrightarrows\mathbb{R}^n$ be a convex process. We define \textit{negative} and \textit{positive dual} processes $H^-$ and $H^+$ of $H$ as follows:
\begin{subequations}
	\begin{align} \label{eq:def of dual} 
	p\in H^-(q) &\iff \langle p,x\rangle\geq \langle q,y\rangle\quad \forall\, (x,y)\in\gr(H),  \\
	p\in H^+(q) & \iff \langle p,x\rangle\leq \langle q,y\rangle\quad \forall\, (x,y)\in\gr(H). 
	\end{align}
\end{subequations}

Note that $H^+(q)= -H^-(-q)$ for all $q$. If $H$ is a closed convex process, we know that $(H^+)^- = H$ and
\begin{equation}\label{eq:H(0) to domain} H(0)=(\dom H^+)^+=(\dom H^-)^-.\end{equation} 
It is straightforward to check that 
\begin{equation}\label{eq:dual in graphs} \graph(H^-) \!=\! \begin{bmatrix} 0 &\!\! I_n \\ -I_n &\!\! 0 \end{bmatrix} \left(\graph H\right)^-\!, \graph(H^+)\!=\! \begin{bmatrix} 0 &\!\! I_n \\ -I_n &\!\! 0 \end{bmatrix} \left(\graph H\right)^+\!\!. \end{equation} 
For a convex cone $ \calC\subseteq\R^n$, we define $\lin(\calC) =-\calC\cap \calC$ and $\Lin(\calC)=\calC-\calC$. We can now define two linear processes $L_-$ and $L_+$ associated with $H$ by
\begin{equation}\label{eq:def of L- and L+}
\graph(L_-):=\lin\big(\gr(H)\big)\textrm{ and } \graph(L_+):=\Lin\big(\graph(H)\big).
\end{equation}
By definition, we therefore have
\begin{equation}\label{eq:lh-+}
\graph(L_-)\subseteq\graph(H)\subseteq\graph(L_+).
\end{equation}
It is clear that $L_-$ and $L_+$ are, respectively, the largest and the smallest (with respect to the graph inclusion) linear processes satisfying \eqref{eq:lh-+}. We call $L_-$ and $L_+$, respectively, the minimal and maximal linear processes associated with $H$. If $H$ is not clear from context, we write $L_-(H)$ and $L_+(H)$ in order to avoid confusion. 

If $L$ is a linear process it is clear that the negative and positive dual processes are equal, which allows us to denote it by $L^\bot:=L^- =L^+$. In fact, the minimal and maximal linear processes associated with a convex process enjoy the following additional properties:
\begin{subequations}\label{eq:L(H) vs L(H+-)}
	\begin{align}
	&L_-(H^-)=L_-(H^+)=L_+^\bot, \\
	&L_+(H^-)=L_+(H^+)=L_-^\bot.
	\end{align}
\end{subequations}
For the reachable and null-controllable sets of $L_-$ and $L_+$ we use the following shorthand notation: 
\begin{align*}
\calR_- := \calR(L_-), \quad &\calR_+:= \calR(L_+). \\
\calN_- := \calN(L_-), \quad &\calN_+:= \calN(L_+). 
\end{align*}

We denote the image of a set $\calS$ under a convex process $H$ by $H(\calS) := \{ y\in\mathbb{R}^n \mid \exists x\in \calS \textrm{ s.t. } y\in H(x)\}$. A direct consequence of this definition is that
\begin{equation}\label{eq:H(S) in graph} H(\calS) = \begin{bmatrix} 0 & I_n \end{bmatrix} \big( \gr(H) \cap \left(\calS\times\mathbb{R}^n\right)\big). \end{equation}
We can define powers of convex processes, by taking $H^0$ equal to the identity map, and letting for $q\geq 0$: 
\[H^{q+1}(x) := H(H^{q}(x)) \quad \forall x\in \mathbb{R}^n. \]
We can define the inverse of a convex process by $H^{-1}(y)=\{ x \mid y\in H(x)\}$. Note that this is always defined as a set-valued map. For higher negative powers of $H$ we use the shorthand: $H^{-n}(x) = (H^{-1})^n(x)$.

Let $L:\mathbb{R}^n\rightrightarrows\mathbb{R}^n$ be a linear process, then we know that $\calF(L)= L^{-n}(\mathbb{R}^n)$ and $\calR(L)= L^n(0)$. In addition
\begin{subequations}\label{eq:LorthRF}
	\begin{gather} 
	\calF(L^\bot) = \calR(L)^\bot, \\
	\calR(L^\bot) = \calF(L)^\bot.
	\end{gather}
\end{subequations}

%
We will characterize reachability in terms of spectral conditions. For this we require one more definition: A real number $\lambda$ and vector $\xi\in\R^n\setminus\{0\} $ form an \textit{eigenpair} of $H$ if $\lambda\xi\in H(\xi)$. In this case $\lambda$ is called an \textit{eigenvalue} and $\xi$ is called an \textit{eigenvector} of $H$. 

In the following, we will need the assumption:
\begin{equation}\label{eq:reach assump} \dom H +\calR_- =\R^n. \end{equation}
As proven in \cite[Thm. 1, Lem. 7]{ReachNullc:19}, we can characterize reachability in terms of eigenvalues of the dual process.

\begin{theorem}\label{thm:cp reach}
	Let $H$ be a convex process such that \eqref{eq:reach assump} holds. Then, the following are equivalent:
	\begin{enumerate}
		\item $H$ is reachable.
		\item $\calR(H)=\R^n$.
		\item $\calR_+=\R^n$ and $H^{-}$ has no nonnegative eigenvalues.
	\end{enumerate}
\end{theorem}
We now move towards null-controllability. It is tempting to think that null-controllability of $H$ is equivalent to reachability of $H^{-1}$. However, while indeed it is true that $\calR(H^{-1})=\calN(H)$, we do not necessarily have that $\calF(H^{-1}) = \calF(H)$. 

As such, we require a characterization of null-controllability. This will be done under slightly more restrictive assumptions than Theorem~\ref{thm:cp reach}. To be precise, we will assume both \eqref{eq:reach assump} and
\begin{equation}\label{eq:nullc assump}
	\rrp = \im H+\calN_-=\R^n.
\end{equation}
The following was proven in \cite[Thm. 2, Lem. 9]{ReachNullc:19}: 
\begin{theorem}\label{thm:nullc}
	Let $H$ be a convex process such that \eqref{eq:reach assump} and \eqref{eq:nullc assump} hold. Then, the following are equivalent:
	\begin{enumerate}
		\item $H$ is null-controllable.
		\item $\calN(H)-\calR(H)=\mathbb{R}^n$. 
		\item $H^{-}$ has no positive eigenvalues. 
	\end{enumerate}
\end{theorem}

The following shows why we require separate tests for these two properties. 
\begin{example} 
	Recall that, as is the case for discrete time linear systems, a convex process can be null-controllable without being reachable. As a simple example consider the convex process given by:
	\[ \graph H := \mathbb{R}\times \{0\}. \] 
	On the other hand, we know that reachability implies null-controllability for discrete time linear systems. For general convex processes this is not the case. As an example, let:
	\[ \graph G := \{ (x,y) \mid 0\leq x\leq y \}. \]
	Note that $\calR(G)=\mathbb{R}_+=\calF(G)$, and therefore $G$ is reachable. As any trajectory of $G$ is a non-decreasing sequence, $G$ is clearly not null-controllable. This means that in general tests for reachability can not be employed to obtain results for null-controllability. 
\end{example}

These two theorems allow us to check for reachability and null-controllability without explicitly determining $\calR(H)$ or $\calN(H)$. This will be central in resolving Problem~\ref{prob: prob} in the next section.  
\vspace{-0.5em}\section{Informativity for convex processes}\label{sec:inf}
We turn our attention to the context of informativity. Let $\calD\subseteq\mathbb{R}^n\times \mathbb{R}^n$ be a finite set of measurements. We define the \textit{most powerful unfalsified process}, $H_\calD$, by:
\[ \graph H_\calD := \cone \calD. \]
By definition we see that $H_\calD\in\Sigma_\calD$ and $\graph H_\calD \subseteq \graph H$ if and only if $H\in\Sigma_\calD$. Our goal is to find conditions on $\calD$ under which every $H\in\Sigma_\calD$ is reachable or null-controllable. we start with the following theorem:

\begin{theorem}\label{thm:reach iff Rn}
	Suppose that \eqref{eq:reach assump} holds for $H_\calD$. Then $H_\calD$ is reachable if and only if every $H\in\Sigma_\calD$ is reachable. 
\end{theorem}
\BP Note that $H_\calD\in\Sigma_\calD$. Therefore the `if' part is immediate. For the `only if' part, assume that $H_\calD$ is reachable. By Theorem~\ref{thm:cp reach}, we have that $\calR(H_\calD)=\mathbb{R}^n$. Now let $H$ be a convex process such that $\graph H_\calD \subseteq \graph H$. As any $q-$step trajectory of $H_\calD$ is one of $H$, it is immediate that $\calR(H_\calD)\subseteq\calR(H)$. Therefore $\calR(H)=\mathbb{R}^n$. This implies that $H$ is reachable. \EP
\vspace{-1em}\begin{remark}\label{rem:r not rn}
It is important to stress that a convex process $H$ is defined to be reachable if $\calF(H)\subseteq\calR(H)$. Therefore a nonstrict convex process $H$ can be reachable whilst $\calR(H)\neq \mathbb{R}^n$. Now let $\graph H\subseteq \graph G$. Note that we may \textit{not} conclude reachability of $G$ from reachability of $H$ in general. As an example, let $\graph H = \{0\}$. This convex process is reachable, and its graph is contained in the graph of any other convex process, which are not necessarily reachable.
\end{remark}

Next, we study null-controllability. It is clear that the reasoning of Remark~\ref{rem:r not rn} also applies to null-controllability. This leads to an important point of contrast between Theorem~\ref{thm:cp reach} and Theorem~\ref{thm:nullc}: Under the conditions of the latter the convex process $H$ can be null-controllable even if $\calN(H)\neq \mathbb{R}^n$. 

\begin{theorem}\label{thm:nullc stab iff Rn}
	Suppose that \eqref{eq:reach assump} and \eqref{eq:nullc assump} hold for $H_\calD$. Then, $H_\calD$ is null-controllable if and only if every $H\in\Sigma_\calD$ is null-controllable. 
\end{theorem}
\BP Again the `if' part is immediate. For the `only if' part, assume that $H_\calD$ is null-controllable. Let $H$ be a convex process such that $\graph H_\calD \subseteq \graph H$. As in the proof of Theorem~\ref{thm:reach iff Rn}, we see that $\calR(H_\calD)\subseteq\calR(H)$ and $\calN(H_\calD)\subseteq\calN(H)$. This implies that 
\[ \mathbb{R}^n= \calN(H_\calD)-\calR(H_\calD) \subseteq \calN(H)-\calR(H).\]

Note that we also have $\graph L_-(H_\calD) \subseteq \graph L_-(H)$ and $\graph L_+(H_\calD) \subseteq \graph L_+(H)$. Therefore, it is clear that \eqref{eq:reach assump} and \eqref{eq:nullc assump} hold for $H$. This implies that $H$ is null-controllable. \EP
\vspace{-0.5em}
The question rests whether we can provide simple tests for reachability and null-controllability of $H_\calD$ in terms of the data $\calD$. In order to resolve this, we will begin by giving two equivalent representations of $H_\calD$. 

Denote $T= |\calD|$ and $\calD = \{(x_t,y_t): t=1,\ldots,T\}$. We define the matrices $X,Y\in\mathbb{R}^{n\times T}$ by taking: 
\[X:= \bbm x_1 & x_2 & \cdots & x_T \ebm,\quad Y := \bbm y_1 & y_2 & \cdots & y_T \ebm.\]
Since $\cone\calD$ is a convex cone, we have that $\calD^+=(\cone \calD)^+$. As $\calD$ is a finite set, we have that $\cone\calD$ and $\calD^+$ are polyhedral cones. This means that there exists $\ell\in\mathbb{N}$ and $\eta_1,\ldots,\eta_{\ell}\in\mathbb{R}^{2n}$, such that $\calD^+=\cone\{\eta_1,\ldots,\eta_{\ell}\}$. We can now define matrices $Z,W\in\mathbb{R}^{\ell\times n}$ by the following partition: 
\[ \begin{bmatrix} Z & -W \end{bmatrix} := \begin{bmatrix} \eta_1  &\ldots & \eta_{\ell} \end{bmatrix}^\top. \]
As $\cone \calD $ is closed, it is equal to $(\calD^+)^+$. Recall that $\graph H_\calD = \cone \calD$. Therefore, we can use \eqref{eq:conic hull rep} to represent $H_\calD$ in the following ways:
\begin{equation}\label{eq:reps of HD} \graph H_\calD \!=\! \begin{bmatrix} X \\ Y \end{bmatrix} \mathbb{R}_+^T=\! \left\lbrace (x,y) \mid \begin{bmatrix} Z & \!-W \end{bmatrix} \begin{bmatrix} x\\y\end{bmatrix} \in \mathbb{R}_+^{\ell} \right\rbrace.\end{equation}

Immediately, we see that 
\[ \dom H_\calD = X\mathbb{R}_+^T \qand  \im H_\calD= Y\mathbb{R}_+^T.\]
Using \eqref{eq:reps of HD} we can express the minimal and maximal linear processes of $H_\calD$ as follows:
\vspace{-0.5em}\begin{align*}
\graph L_-(H_\calD) &= \ker \begin{bmatrix} Z & -W \end{bmatrix}, \\
\graph L_+(H_\calD) &= \im \begin{bmatrix} X \\ Y \end{bmatrix}. 
\end{align*}
For the characterizations of reachability and null-controllability in Theorem~\ref{thm:cp reach} and Theorem~\ref{thm:nullc} respectively, we need the reachable and null-controllable sets of $L_+$ and $L_-$. In order to characterize these in terms of the data $\calD$, we first look at the image of a set under these linear processes. For a given set $\calS\subseteq \mathbb{R}^n$ we can apply \eqref{eq:H(S) in graph} to verify that:
\begin{align*}
L_-(H_\calD)(\calS) &= W^{-1}Z\calS, \\
L_+(H_\calD)(\calS) &= YX^{-1}\calS.
\end{align*}
Recall that for a linear process $L$ the reachable set is \textit{finitely determined} and $\calR(L)=L^n(0)$. Combining the above with some slight abuse of notation, we can write:
\begin{align*}
\calR(L_-(H_\calD)) &= (W^{-1}Z)^n\{0\}, \\
\calR(L_+(H_\calD)) &= (YX^{-1})^n\{0\}.
\end{align*}
This characterizes the reachable sets of $L_-(H_\calD)$ and $L_+(H_\calD)$ using subspace algorithms with at most $n$ steps. Following the same reasoning with negative powers, we obtain that: 
\begin{align*}
\calN(L_-(H_\calD)) &= (Z^{-1}W)^n\{0\},  \\
\calN(L_+(H_\calD)) &= (XY^{-1})^n\{0\}.
\end{align*}

We now shift our focus to the negative dual of $H_\calD$, and show that it can be represented in terms of $X$ and $Y$ or $Z$ and $W$ as well. 

By \eqref{eq:dual in graphs} and the first representation of \eqref{eq:reps of HD} we have that: 
\[ \graph(H^-_\calD) = \begin{bmatrix} 0 & I_n \\ -I_n & 0 \end{bmatrix} \left(\begin{bmatrix} X \\ Y \end{bmatrix} \mathbb{R}_+^T \right)^-. \] 
By \eqref{eq:dual and linear map} this implies that: 
\[ \graph(H^-_\calD) \!=\! \begin{bmatrix} 0 &\! I_n \\ -I_n &\! 0 \end{bmatrix}\!\!\begin{bmatrix} X^\top &\!\! Y^\top \end{bmatrix}^{-1}\!\mathbb{R}_-^T=\!\begin{bmatrix} Y^\top &\!\!-X^\top \end{bmatrix}^{-1}\!\mathbb{R}_-^T. \]
Similarly, we can begin from \eqref{eq:dual in graphs} and the second representation in \eqref{eq:reps of HD} instead. As such, we can conclude that the negative dual of $H_\calD$ satisfies:
\[ \graph(H_\calD^-) \!=\! \begin{bmatrix} W^\top \\Z^\top \end{bmatrix}\mathbb{R}_+^\ell =\! \left\lbrace (x,y) \mid \begin{bmatrix} Y^\top &\! -X^\top \end{bmatrix} \begin{bmatrix} x\\y\end{bmatrix} \in \mathbb{R}_-^T \right\rbrace.\]

Then, we have that $\lambda$ and $\xi$ form an eigenpair of $H_\calD^-$ if and only if $\xi\neq 0$ and $\xi^\top(Y -\lambda X) \leq 0$. 

%

We can now combine the previous discussion with Theorem~\ref{thm:cp reach} and Theorem~\ref{thm:reach iff Rn} to obtain the following characterization of informativity for reachability in terms of data: 
\begin{theorem}\label{thm:inf reach}
	Let $\calD\subseteq\mathbb{R}^n\times \mathbb{R}^n$ be a finite set. Suppose that \[X\mathbb{R}_+^T+(W^{-1}Z)^n\{0\}=\mathbb{R}^n. \]
	Then, $\calD$ is informative for reachability if and only if $(YX^{-1})^n\{0\}=\mathbb{R}^n$ and for all $\lambda\geq 0:$
	\[ \xi^\top(Y -\lambda X) \leq 0 \implies \xi = 0. \] 	
\end{theorem}
\begin{remark}
	Note that $(YX^{-1})^n\{0\}=\mathbb{R}^n$ implies that $\calR(L)=\mathbb{R}^n$ for all \textit{linear processes} $L$ such that $\calD \subseteq \graph L$. That is, all such linear processes are reachable. 
\end{remark}

\begin{example}
	Let $n=2$ and suppose that we measure the following $4-$step trajectory: 
	\[ x_0 \!=\! \begin{bmatrix} 0 \\ 0 \end{bmatrix}\!, x_1\!=\!  \begin{bmatrix} 1 \\ 0 \end{bmatrix}\!, x_2\!=\!  \begin{bmatrix} 0 \\ 1 \end{bmatrix}\!, x_3\!=\!  \begin{bmatrix} 0 \\ -1 \end{bmatrix}\!, x_4\!=\!  \begin{bmatrix} -1 \\ 0 \end{bmatrix}. \]
	If we define $X$ and $Y$ as before, we get
	\[ X= \begin{bmatrix} 0 & 1 & 0 & 0 \\ 0 & 0 & 1 & -1\end{bmatrix}, \quad Y= \begin{bmatrix} 1 & 0 & 0 &-1 \\  0 & 1 & -1&0\end{bmatrix}.\]	
	We can use these to find $Z$ and $W$:
	\[\begin{bmatrix} Z & -W \end{bmatrix} = \begin{bmatrix} 1 & 0 & 0 & 0 \\ 1 & 0 & 0 & -1\end{bmatrix}. \]
	First, note that $X\mathbb{R}^4_+ = \mathbb{R}_+ \times \mathbb{R}$ and $(W^{-1}Z)^2\{0\} = \mathbb{R}\times \{0\}$. Therefore, we can now use Theorem~\ref{thm:inf reach} to check for informativity. 
	
	Now, it is straightforward to verify that $(YX^{-1})^2\{0\}=\mathbb{R}^2$. Lastly, let $\lambda \geq 0$ and  \[\begin{bmatrix} \xi_1 & \xi_2 \end{bmatrix}\begin{bmatrix}  1 & -\lambda  & 0 & -1 \\ 0 & 1 & -1-\lambda  & \lambda \end{bmatrix} \leq 0.\]
	By direct inspection, it is clear that this implies that 
	\[ \xi_1\leq 0, \quad \xi_2 \leq \lambda \xi_1, \quad 0 \leq (1+\lambda)\xi_2, \quad \lambda \xi_2 \leq \xi_1.\]
	These inequalities show that for any $\lambda \geq 0$ we have that $\xi_1=\xi_2=0$. This proves that $\calD$ is informative for reachability. 
\end{example}

In a similar fashion we can apply our discussion to Theorem~\ref{thm:nullc} and Theorem~\ref{thm:nullc stab iff Rn} to obtain a characterization of informativity for null-controllability. 

\begin{theorem}\label{thm:inf null-c}
	Let $\calD\subseteq\mathbb{R}^n\times \mathbb{R}^n$ be a finite set. Suppose that
\[X\mathbb{R}_+^T+(W^{-1}Z)^n\{0\} =\mathbb{R}^n \]
and
\[(YX^{-1})^n\{0\}=Y\mathbb{R}_+^T+(Z^{-1}W)^n\{0\}= \mathbb{R}^n.\]
Then $\calD$ is informative for null-controllability if and only if for all $\lambda> 0:$
\[ \xi^\top(Y -\lambda X) \leq 0 \implies \xi = 0. \] 
\end{theorem}

\begin{remark} 
	If $H$ is a convex process whose graph is polyhedral, we can always find a finite set $\calD$ such that $H=H_\calD$. This means that the results of Theorem~\ref{thm:inf reach} and Theorem~\ref{thm:inf null-c} can be applied to any polyhedral convex process without loss of generality.
\end{remark}
\vspace{-0.5em}\section{Conclusions}\label{sec:conc}
In this paper, we have resolved a number of informativity problems for conically constrained linear systems. This means that we have formulated conditions on finite, exact, state measurements under which we can test whether the measured system is reachable or null-controllable. The resulting tests take the convenient form of subspace inclusions and spectral conditions. 

Future work includes extending the ideas in this paper towards the more general class of linear systems with convex constraints. It is easy to see that these systems can be viewed as difference inclusions of convex set-valued maps. Similar to the approach in this paper, we can define the smallest set-valued map consistent with the data by taking the convex hull instead of the conic hull. As such, a characterization of reachability for such systems will lead to informativity results for this class of systems. Another direction of future work is investigating informativity for the analysis of other properties and for control. Interesting problems are for example dissipativity or feedback stabilization. Resolving such a problem would require formulating characterizations for a given convex process to have the aforementioned properties. Lastly, this paper considers only exact measurements of the state. However, many realistic scenarios will involve noisy measurements. Incorporating noisy data within this framework will lead to interesting informativity problems. 
	
\bibliographystyle{IEEEtran}
\bibliography{references}

\end{document}

%% file: ka-newcommands.tex
\DeclareMathOperator{\im}{im}

\DeclareMathOperator{\cl}{cl}

\DeclareMathOperator{\cone}{cone}

\DeclareMathOperator{\dom}{dom}
\DeclareMathOperator{\graph}{gr}
\DeclareMathOperator{\gr}{gr}
\DeclareMathOperator{\lin}{lin}
\DeclareMathOperator{\Lin}{Lin}

\let\leq\leqslant
\let\geq\geqslant
\let\emptyset\varnothing

\newcommand{\calC}{\ensuremath{\mathcal{C}}}
\newcommand{\calD}{\ensuremath{\mathcal{D}}}

\newcommand{\calF}{\ensuremath{\mathcal{F}}}

\newcommand{\calN}{\ensuremath{\mathcal{N}}}

\newcommand{\calR}{\ensuremath{\mathcal{R}}}
\newcommand{\calS}{\ensuremath{\mathcal{S}}}
\newcommand{\calT}{\ensuremath{\mathcal{T}}}

\newcommand{\bmat}{\begin{matrix}}
\newcommand{\emat}{\end{matrix}}
\newcommand{\bbm}{\begin{bmatrix}}
\newcommand{\ebm}{\end{bmatrix}}
\newcommand{\bbma}{\begin{bmatrix*}[r]}
\newcommand{\ebma}{\end{bmatrix*}}
\newcommand{\bpm}{\begin{pmatrix}}
\newcommand{\epm}{\end{pmatrix}}
\newcommand{\bvm}{\begin{vmatrix}}
\newcommand{\evm}{\end{vmatrix}}
\newcommand{\bse}{\begin{subequations}}
\newcommand{\ese}{\end{subequations}}
\newcommand{\beq}{\begin{equation}}
\newcommand{\eeq}{\end{equation}}
\newcommand{\ben}{\renewcommand{\labelenumi}{\arabic{enumi}.}
\renewcommand{\theenumi}{\arabic{enumi}}\begin{enumerate}}

\newcommand{\een}{\end{enumerate}}

\newcommand{\beni}{\renewcommand{\labelenumi}{\roman{enumi}.}
\renewcommand{\theenumi}{\roman{enumi}}\begin{enumerate}}

\newcommand{\eeni}{\end{enumerate}}

\newcommand{\bena}{\renewcommand{\labelenumi}{\alph{enumi}.}
\renewcommand{\theenumi}{\alph{enumi}}\begin{enumerate}}

\newcommand{\eena}{\end{enumerate}}

\newcommand{\bit}{\begin{itemize}}
\newcommand{\eit}{\end{itemize}}
\newcommand{\bthe}{\begin{theorem}}
\newcommand{\ethe}{\end{theorem}}
\newcommand{\blem}{\begin{lemma}}
\newcommand{\elem}{\end{lemma}}
\newcommand{\bprop}{\begin{proposition}}
\newcommand{\eprop}{\end{proposition}}
\newcommand{\bex}{\begin{example}}
\newcommand{\eex}{\end{example}}
\newcommand{\bas}{\begin{assumption}}
\newcommand{\eas}{\end{assumption}}
\newcommand{\bre}{\begin{remark}}
\newcommand{\ere}{\end{remark}}
\newcommand{\bcor}{\begin{corollary}}
\newcommand{\ecor}{\end{corollary}}
\newcommand{\bdfn}{\begin{definition}}
\newcommand{\edfn}{\end{definition}}
\newcommand{\bcon}{\begin{conjecture}}
\newcommand{\econ}{\end{conjecture}}

\newcommand{\R}{\ensuremath{\mathbb R}}

\newcommand{\N}{\ensuremath{\mathbb N}}
\newcommand{\Z}{\ensuremath{\mathbb Z}}

\newcommand{\BP}{\noindent{\bf Proof. }}
\newcommand{\EP}{\hspace*{\fill} $\blacksquare$\bigskip\noindent}
\newcommand{\qand}{\quad\text{ and }\quad}

%% file: root.bbl
\begin{thebibliography}{10}
\providecommand{\url}[1]{#1}
\csname url@samestyle\endcsname
\providecommand{\newblock}{\relax}
\providecommand{\bibinfo}[2]{#2}
\providecommand{\BIBentrySTDinterwordspacing}{\spaceskip=0pt\relax}
\providecommand{\BIBentryALTinterwordstretchfactor}{4}
\providecommand{\BIBentryALTinterwordspacing}{\spaceskip=\fontdimen2\font plus
\BIBentryALTinterwordstretchfactor\fontdimen3\font minus
  \fontdimen4\font\relax}
\providecommand{\BIBforeignlanguage}[2]{{%
\expandafter\ifx\csname l@#1\endcsname\relax
\typeout{** WARNING: IEEEtran.bst: No hyphenation pattern has been}%
\typeout{** loaded for the language `#1'. Using the pattern for}%
\typeout{** the default language instead.}%
\else
\language=\csname l@#1\endcsname
\fi
#2}}
\providecommand{\BIBdecl}{\relax}
\BIBdecl

\bibitem{vanWaarde2020}
H.~J. {van Waarde}, J.~{Eising}, H.~L. {Trentelman}, and M.~K. {Camlibel},
  ``Data informativity: a new perspective on data-driven analysis and
  control,'' \emph{IEEE Transactions on Automatic Control}, vol.~65, no.~11,
  pp. 4753--4768, 2020.

\bibitem{Eising2020c}
J.~Eising and H.~L. Trentelman, ``Informativity of noisy data for structural
  properties of linear systems,'' \emph{https://arxiv.org/abs/2009.01552},
  2020.

\bibitem{vanWaarde2020c}
H.~J. van Waarde, M.~K. Camlibel, and M.~Mesbahi, ``From noisy data to feedback
  controllers: non-conservative design via a matrix {S}-lemma,'' \emph{IEEE
  Transactions on Automatic Control, early access}, 2020.

\bibitem{Trentelman2020}
H.~L. Trentelman, H.~J. {van Waarde}, and M.~K.~. {Camlibel}, ``An
  informativity approach to data-driven tracking and regulation,''
  \emph{https://arxiv.org/abs/2009.01552}, 2020.

\bibitem{Koch2020}
A.~{Koch}, J.~{Berberich}, and F.~{Allg{\"o}wer}, ``Verifying dissipativity
  properties from noise-corrupted input-state data,'' in \emph{Proceedings of
  the IEEE Conference on Decision and Control}, 2020, pp. 616--621.

\bibitem{DePersis2020}
C.~{De Persis} and P.~{Tesi}, ``Formulas for data-driven control:
  Stabilization, optimality, and robustness,'' \emph{IEEE Transactions on
  Automatic Control}, vol.~65, no.~3, pp. 909--924, 2020.

\bibitem{Berberich2020}
J.~{Berberich}, A.~{Koch}, C.~W. {Scherer}, and F.~{Allg\"ower}, ``Robust
  data-driven state-feedback design,'' in \emph{Proceedings of the American
  Control Conference}, 2020, pp. 1532--1538.

\bibitem{Angeli:09}
D.~Angeli, P.~De~Leenheer, and E.~D. Sontag, ``Chemical networks with inflows
  and outflows: A positive linear differential inclusions approach,''
  \emph{Biotechnology Progress}, vol.~25, no.~3, pp. 632--642, 2009.

\bibitem{makarov:77}
V.~L. Makarov and A.~M. Rubinov, \emph{Mathematical Theory of Economic Dynamics
  and Equilibria}.\hskip 1em plus 0.5em minus 0.4em\relax Springer-Verlag,
  1977.

\bibitem{j37}
M.~D. Kaba and M.~K. Camlibel, ``A spectral characterization of controllability
  for linear discrete-time systems with conic constraints,'' \emph{SIAM Journal
  on Control and Optimization}, vol.~53, no.~4, pp. 2350--2372, 2015.

\bibitem{oh:05-2}
S.~R. Oh and S.~K. Agrawal, ``A reference governor-based controller for a cable
  robot under input constraints,'' \emph{IEEE Transactions on Control Systems
  Technology}, vol.~13, no.~4, pp. 639--645, 2005.

\bibitem{Frankowska:87}
H.~Frankowska, ``Local controllability and infinitesimal generators of
  semigroups of set-valued maps,'' \emph{SIAM Journal on Control and
  Optimization}, vol.~25, no.~2, pp. 412--432, 1987.

\bibitem{Aubin:84}
J.-P. Aubin and I.~Ekeland, \emph{Applied nonlinear analysis}, ser. Pure and
  applied mathematics.\hskip 1em plus 0.5em minus 0.4em\relax J. Wiley, 1984.

\bibitem{AFO:86}
J.-P. Aubin, H.~Frankowska, and C.~Olech, ``Controllability of convex
  processes,'' \emph{SIAM Journal on Control and Optimization}, vol.~24, no.~6,
  pp. 1192--1211, 1986.

\bibitem{PD:94}
V.~N. Phat and T.~C. Dieu, ``On the {K}re\u\i n-{R}utman theorem and its
  applications to controllability,'' \emph{Proceedings of the American
  Mathematical Society}, vol. 120, no.~2, pp. 495--500, 1994.

\bibitem{ReachNullc:19}
J.~Eising and M.~K. Camlibel, ``On reachability and null-controllability of
  nonstrict convex processes,'' \emph{IEEE Control Systems Letters}, vol.~3,
  no.~3, pp. 751--756, 2019.

\bibitem{Park2009}
U.~S. Park and M.~Ikeda, ``Stability analysis and control design of {LTI}
  discrete-time systems by the direct use of time series data,''
  \emph{Automatica}, vol.~45, no.~5, pp. 1265--1271, 2009.

\bibitem{Wang2011}
Z.~{Wang} and D.~{Liu}, ``Data-based controllability and observability analysis
  of linear discrete-time systems,'' \emph{IEEE Transactions on Neural
  Networks}, vol.~22, no.~12, pp. 2388--2392, Dec 2011.

\bibitem{Liu2014}
D.~{Liu}, P.~{Yan}, and Q.~{Wei}, ``Data-based analysis of discrete-time linear
  systems in noisy environment: Controllability and observability,''
  \emph{Information Sciences}, vol. 288, pp. 314--329, 12 2014.

\bibitem{Niu2017}
H.~{Niu}, H.~{Gao}, and Z.~{Wang}, ``A data-driven controllability measure for
  linear discrete-time systems,'' in \emph{IEEE 6th Data Driven Control and
  Learning Systems Conference}, May 2017, pp. 455--466.

\bibitem{Zhou2018}
B.~{Zhou}, Z.~{Wang}, Y.~{Zhai}, and H.~{Yuan}, ``Data-driven analysis methods
  for controllability and observability of a class of discrete {LTI} systems
  with delays,'' in \emph{IEEE 7th Data Driven Control and Learning Systems
  Conference}, May 2018, pp. 380--384.

\bibitem{Mayne:00}
D.~Q. Mayne, J.~B. Rawlings, C.~V. Rao, and P.~O.~M. Scokaert, ``Constrained
  model predictive control: Stability and optimality,'' \emph{Automatica},
  vol.~36, pp. 789--814, 2000.

\bibitem{Mayne2014}
D.~Q. Mayne, ``Model predictive control: Recent developments and future
  promise,'' \emph{Automatica}, vol.~50, no.~12, pp. 2967 -- 2986, 2014.

\bibitem{Piga2018}
D.~{Piga}, S.~{Formentin}, and A.~{Bemporad}, ``Direct data-driven control of
  constrained systems,'' \emph{IEEE Transactions on Control Systems
  Technology}, vol.~26, no.~4, pp. 1422 -- 1429, 2018.

\bibitem{Coulson2019}
J.~{Coulson}, J.~{Lygeros}, and F.~{D\"orfler}, ``Data-enabled predictive
  control: In the shallows of the {DeePC},'' in \emph{Proceedings of the
  European Control Conference}, June 2019, pp. 307--312.

\bibitem{Aubin:90}
J.-P. Aubin and H.~Frankowska, \emph{Set-valued Analysis}, ser. Systems \&
  Control: Foundations \& Applications.\hskip 1em plus 0.5em minus 0.4em\relax
  Boston, MA: Birkh\"auser Boston Inc., 1990, vol.~2.

\end{thebibliography}
